\documentclass[10pt]{article}

\usepackage{fullpage}
\usepackage{amsmath}
\usepackage{amssymb}
\usepackage{dsfont}
\usepackage{amsfonts}
\usepackage{hyperref}
\usepackage{graphicx}
\usepackage{epsfig}
\usepackage[authoryear,round]{natbib}
\usepackage{array}

\newtheorem{remark}{Remark} 

\date{}
\begin{document}

\title{Global Sensitivity Analysis with Dependence Measures}

\author{Sebastien Da Veiga$^{a}$\thanks{Now at Snecma (Safran), \'Etablisssement de Montereau, A\'erodrome Melun-Villaroche, BP 1936, 77019 Melun cedex, France. Email: sebastien.daveiga@snecma.fr \vspace{6pt}}
\\ \vspace{6pt}  
$^{a}${\em{IFP Energies nouvelles,
1 \& 4 avenue de Bois Pr\'{e}au, 92852 Rueil-Malmaison, France}}}

\maketitle

\begin{abstract}
Global sensitivity analysis with variance-based measures suffers from several theoretical and practical limitations,  
since they focus only on the variance of the output and handle multivariate variables in a limited way. In this paper, we introduce 
a new class of sensitivity indices based on dependence measures which overcomes these insufficiencies. Our approach originates 
from the idea to compare the output distribution with its conditional counterpart when one of the input variables is fixed. We establish that 
this comparison yields previously proposed indices when it is performed with Csisz\'ar f-divergences, as well as sensitivity indices which are well-known dependence measures between random variables. This leads us to investigate completely new sensitivity indices based on recent state-of-the-art dependence measures, such as distance correlation and the Hilbert-Schmidt independence criterion. We also emphasize the potential of feature selection techniques relying on such dependence measures as alternatives to screening in high dimension.
\end{abstract}

\section{Introduction}

Since the early work of Sobol \citep{Sobol93}, global sensitivity analysis (GSA)  has received a lot of attention in the computer code 
experiments community. These days variance-based sensitivity indices are common tools in the analysis of complex 
physical phenomena. Several statistical estimators have been proposed \citep{cukier73,saltelli10, tarantola06, janon13, owen13} and their asymptotic 
properties are now well understood \citep{tissot12, sdv13, janon13}. In addition, the case of 
computationally expensive codes has been investigated thoroughly with the introduction of several dedicated surrogate models 
\citep{OOH04, marrel09, sdv09, blatman10, touzani12, durrande12}.

However, even if they are extremely popular and informative importance measures, variance-based indices suffer from theoretical and practical limitations. First, by definition they only study the impact of the input parameters on the variance of the output. Since this is a restricted summary of the output distribution, this measure happens to be inadequate for many case studies. Alternative approaches include for example density-based indices \citep{borgonovo07}, derivative-based measures \citep{sobol09}, or goal-oriented dedicated indices \citep{fort13}.
Second, variance-based indices do not generalize easily to the case of a multivariate output \citep{gamboa13}. Unfortunately, computer code outputs often consist of several scalars or even time-dependent curves, which limits severely the practical use of standard indices.
Finally for high-dimensional problems, a preliminary screening procedure is usually mandatory before the analysis of the computer code or the modeling with a surrogate. The computational cost of GSA is in general too high to envision its use for screening purposes and more qualitative approaches are thus needed, e.g. the Morris method \citep{morris91} or group screening \citep{moon12}.

In this paper, we propose a completely original point of view in order to overcome these limitations. Starting from the general 
framework of GSA and the concept of dissimilarity measure, we introduce a new class of sensitivity indices which comprises as a special case the density-based index of \citet{borgonovo07}. We propose an estimation procedure relying on density ratio estimation and show that it gives access to several different indices for the same computational cost. More importantly, we highlight that other special cases lead to well-known dependence measures, including the mutual information. This link motivates us to investigate the potential of recent state-of-the-art dependence measures as new sensitivity indices, such as the distance correlation \citep{dcor07} or the Hilbert-Schmidt independence criterion \citep{gretton05a}. An appealing property of such measures is that they can handle multivariate random variables very easily. We also discuss how feature selection methods based on these measures can effectively replace standard screening procedures.

The structure of the paper is as follows. In Section \ref{dissim}, we introduce the general GSA framework based on dissimilarity measures and discuss the use of Csisz\'ar f-divergences. We also emphasize the link with mutual information and propose an estimation procedure. In Section \ref{Indep}, we give a review of some dependence measures and show how they can be used as new sensitivity indices. We also provide examples of feature selection techniques in which they are involved. Screening will then be seen as an equivalent to feature selection in machine learning. Finally, several numerical experiments are conducted in Section \ref{secexp} on both analytical and industrial applications. In particular, we illustrate the potential of dependence measures for GSA.

\section{From dissimilarity measures to sensitivity indices} \label{dissim}

Denote $Y = \eta(X^1,\ldots,X^p)$ the computer code output which is a function of the $p$ input random variables $X^k$, $k=1,\ldots,p$ 
where $\eta:\mathbb{R}^p\rightarrow \mathbb{R}$ is assumed to be continuous.
In standard global sensitivity analysis, it is further assumed that the $X^k$ have a known distribution and are independent.
As pointed out by \citet{baucells13}, a natural way of defining the impact of a given input $X^k$ on $Y$ is to consider a 
function which measures the similarity between the distribution of $Y$ and that of $Y|X^k$. 
More precisely, the impact of $X^k$ on $Y$ is given by
\begin{equation}
S_{X^k} = \mathbb{E}_{X^k} \left(d(Y,Y|X^k)\right) \label{Si}
\end{equation}
where $d(\cdot,\cdot)$ denotes a dissimilarity measure between two random variables. 
The advantage of such a formulation is that many choices for $d$ are available, and we will see in what follows that some natural 
dissimilarity measures yield sensitivity indices related to well known quantities.
However before going further, let us note that the naive dissimilarity measure 
\begin{equation}
d(Y,Y|X^k) = \left(\mathbb{E}(Y) - \mathbb{E}(Y|X^k)\right)^2 \label{naive}
\end{equation}
where random variables are compared only through their mean values produces the unnormalized Sobol first-order sensivity index 
$S_{X^k}^1 = \textrm{Var}\left(\mathbb{E}(Y|X^k)\right)$.

\subsection{Csisz\'ar f-divergences} \label{fdiv}
Assuming all input random variables have an absolutely continuous distribution with respect to the Lebesgue measure on $\mathbb{R}$,
 the f-divergence \citep{fdiv67} between $Y$ and $Y|X^k$ is given by
\begin{equation*}
d_f(Y||Y|X^k) = \int_{\mathbb{R}} f\left(\frac{p_Y(y)}{p_{Y|X^k}(y)}\right)p_{Y|X^k}(y) dy
\end{equation*}
where $f$ is a convex function such that $f(1)=0$ and $p_Y$ and $p_{Y|X^k}$ are the probability distribution functions of $Y$ and $Y|X^k$, respectively.
Standard choices for function $f$ include for example
\begin{itemize}
\item Kullback-Leibler divergence: $f(t) = -\ln(t)$ or $f(t) = t\ln(t)$;
\item Hellinger distance: $f(t) = \left(\sqrt{t}-1\right)^2$;
\item Total variation distance: $f(t) = \vert t - 1 \vert$;
\item Pearson $\chi^2$ divergence: $f(t) = (t-1)^2$ or $f(t)= t^2-1$;
\item Neyman $\chi^2$ divergence: $f(t) = (t-1)^2/t$ or $f(t) = (1-t^2)/t$.\\
\end{itemize}

Plugging this dissimilarity measure in (\ref{Si}) yields the following sensitivity index:
\begin{equation}
S_{X^k}^f = \int_{\mathbb{R}^2} f\left(\frac{p_Y(y)p_{X^k}(x)}{p_{X^k,Y}(x,y)}\right)p_{X^k,Y}(x,y) dxdy \label{Sif}
\end{equation}
where $p_{X^k}$ and $p_{X^k,Y}$ are the probability distribution functions of $X^k$ and $(X^k,Y)$, respectively.
First of all, note that inequalities on Csisz\'ar f-divergences imply that such sensitivity indices are positive and equal zero when $Y$ and $X_k$ are independent. 
Also, it is important to note that given the form of $S_{X^k}^f$, it is invariant under any smooth and uniquely invertible 
transformation of the variables $X^k$ and $Y$, see the proof for mutual information in \citet{kraskov04}. This is a major advantage 
over variance-based Sobol sensitivity indices, which are only invariant under linear transformations.\\

It is easy to see that the total variation distance with $f(t) = \vert t - 1 \vert$ gives a sensivity index equal to the one proposed by \citet{borgonovo07}:
\begin{equation*}
S_{X^k} ^f= \int_{\mathbb{R}^2} \lvert p_Y(y)p_{X^k}(x) - p_{X^k,Y}(x,y)\rvert dxdy.
\end{equation*}
In addition, the Kullback-Leibler divergence with $f(t) = -\ln(t)$ yields
\begin{equation*}
S_{X^k}^f = \int_{\mathbb{R}^2} p_{X^k,Y}(x,y) \ln\left(\frac{p_{X^k,Y}(x,y)}{p_Y(y)p_{X^k}(x)}\right)dxdy,
\end{equation*}
that is the mutual information $I(X^k;Y)$ between $X^k$ and $Y$. A normalized version of this sensitivity index was studied by \citet{krzy01}. 
Similarly, the Neyman $\chi^2$ divergence with $f(t) = (1-t^2)/t$ leads to 
\begin{equation*}
S_{X^k}^f = \int_{\mathbb{R}^2} \left(\frac{p_{X^k,Y}(x,y)}{p_Y(y)p_{X^k}(x)}-1\right)^2 p_Y(y)p_{X^k}(x)dxdy,
\end{equation*}
which is the so-called squared-loss mutual information between $X^k$ and $Y$ (or mean square contingency).
These results show that some previously proposed sensitivity indices are actually special cases of more general indices defined through 
Csisz\'ar f-divergences. To the best of our knowledge, this is the first work in which this link is highlighted. Moreover, 
the specific structure of equation (\ref{Sif}) makes it possible to envision more efficient tools for the estimation of these sensitivity indices.
Indeed, it only involves approximating a density ratio rather than full densities. This point is investigated in the next subsection. 
But more importantly, we see that special choices for $f$ define sensivity indices that are actually well-known dependence measures 
such as the mutual information. This paves the way for completely new sensitivity indices based on recent state-of-the-art dependence measures, see Section \ref{Indep}.  

\subsection{Estimation}
Coming back to equation (\ref{Sif}), the goal is to estimate 
\begin{equation*}
S_{X^k}^f = \int_{\mathbb{R}^2} f\left(\frac{1}{r(x,y)}\right)p_{X^k,Y}(x,y) dxdy = \mathbb{E}_{(X^k,Y)}f\left(\frac{1}{r(X^k,Y)}\right)
\end{equation*}
where $r(x,y) = p_{X^k,Y}(x,y)/(p_Y(y)p_{X^k}(x))$ is the ratio between the joint density of $(X^k,Y)$ and the marginals.
Of course, straightforward estimation is possible if one estimates the densities $p_{X^k,Y}(x,y)$, $p_{X^k}(x)$ and $p_Y(y)$ with e.g. 
kernel density estimators. However, it is well known that density estimation suffers from the curse of dimensionality. This limits the 
possible multivariate extensions we discuss in the next subsection. Besides, since only the ratio function $r(x,y)$ is needed, we expect 
more robust estimates by focusing only on it.

Let us assume now that we have a sample $(X^k_i,Y_i)_{i=1,\ldots,n}$ of $(X^k,Y)$, the idea is to build first an estimate $\hat{r}(x,y)$ of the ratio. 
The final estimator $\hat{S}_{X^k}^f$ of $S_{X^k}^f$ will then be given by
\begin{equation}
\hat{S}_{X^k}^f = \frac{1}{n}\sum_{i=1}^n f\left(\frac{1}{\hat{r}(X^k_i,Y_i)}\right).
\end{equation}
Powerful estimating methods for ratios include e.g. maximum-likelihood estimation \citep{sugi08}, unconstrained least-squares importance fitting \citep{sugi09}, among others \citep[see][]{sugibook}. A k-nearest neighbors strategy dedicated to mutual information is also discussed in \cite{kraskov04}.

\subsection{Multivariate extensions}
Given our approach focusing only on densities, it is straightforward to extend the definition of the sensitivity index in equation (\ref{Sif}) 
to any number of input and output variables. We can then study the impact of a given group of input variables $X^{u}=\left\{X^k\right\}_{k\in u}$ 
where $u$ is a subset of $\{1,\ldots,p\}$ on a multivariate output $Y\in \mathbb{R}^q$ with the sensitivity index given by
\begin{equation*}
S_{X^{u}}^f = \int_{\mathbb{R}^{\vert u\vert} \times \mathbb{R}^q} f\left(\frac{p_{Y}(y)p_{X^u}(x)}
{p_{X^u,Y}(x,y)}\right)p_{X^u,Y}(x,y) dxdy.
\end{equation*}
This favorable generalization was already mentioned for the special cases of the total-variation distance and mutual information by \citet{borgonovo07} and \citet{auder08}, respectively.
However, in the high-dimensional setting, estimation of such sensitivity indices is infeasible since even the ratio trick detailed above fails. 
This is thus a severe limitation for screening purposes. We examine efficient alternatives in Section \ref{Indep}.

Moreover, note that extending the naive dissimilarity measure (\ref{naive}) to the multivariate output case naturally leads to consider $ d(Y,Y|X^k) = \Vert\mathbb{E}(Y) - \mathbb{E}(Y|X^k)\rVert_2^2$. Straightforward calculations reveal that the corresponding sensitivity index is then the sum of Sobol first-order sensitivity indices on each output. \citet{gamboa13} showed that this multivariate index is the only one possessing desired invariance properties in the variance-based index family.

\subsection{On the use of other dissimilarity measures}

We focused above on Csisz\'ar f-divergences but other dissimilarity measures exist to compare probability distributions. 
In particular, integral probability metrics \citep[IPM,][]{muller97} are a popular family of distance measures on probabilities given by
\begin{equation}
\gamma_{\mathcal{F}}(\mathbb{P},\mathbb{Q}) = \sup_{f\in\mathcal{F}} \left\vert \int_S fd\mathbb{P} - \int_S fd\mathbb{Q} \right\vert \label{IPM}
\end{equation}
for two probability measures $\mathbb{P}$ and $\mathbb{Q}$ and where $\mathcal{F}$ is a class of real-valued bounded measurable functions on $S$. Just as the choice of function $f$ in Csisz\'ar f-divergences gives rise to different measures, the choice of $\mathcal{F}$ generates different IPMs, e.g. the Wasserstein distance, the Dudley metric or the total variation distance. It is interesting to note that Csisz\'ar f-divergences and IPMs are very distinct classes of measures, since they only intersect at the total variation distance \citep{sri12}. Unfortunately, plugging the general expression (\ref{IPM}) of an IPM in equation (\ref{Si}) no longer yields a closed-form expression for a sensitivity index. However, we plan to study such indices in a future work since estimation of IPMs appears to be easier than for Csisz\'ar f-divergences and is independent of the dimensionality of the random variables \citep{sri12}.

Finally, let us mention the recent work of \citet{fort13} on goal-oriented measures, where they introduce a new class of sensitivity indices
\begin{equation}
S_{X^k}^{\psi} = \mathbb{E}\psi(Y;\theta^*)-\mathbb{E}_{(X^k,Y)}\psi(Y;\theta_k(X^k)) \label{contrast}
\end{equation}
where $\psi$ is the contrast function associated to the features of interest $\theta^* = \arg\min_{\theta} \mathbb{E}\psi(Y;\theta)$ and $\theta_k(x) = \arg\min_{\theta} \mathbb{E}(\psi(Y;\theta)|X^k=x)$ of $Y$ and $Y$ conditionally to $X^k=x$, respectively (note that we only give here the unnormalized version of the index). It is easy to check that (\ref{contrast}) is a special case of (\ref{Si}).

\section{Dependence measures and feature selection} \label{Indep}

Given two random vectors $X$ in $\mathbb{R}^p$ and $Y$ in $\mathbb{R}^q$, dependence measures aim at 
quantifying the dependence between $X$ and $Y$ in arbitrary dimension, with the property that the measure equals zero if 
and only if $X$ and $Y$ are independent. In particular, they are useful when one wants to design a statistical test for independence. 
Here, we focus on the long-known mutual information criterion, as well as on the novel distance correlation measure \citep{dcor07}. 
Recently, \citet{sejdi13} showed that it shares deep links with distances between embeddings of distributions to reproducing kernel Hilbert spaces (RKHS) and especially the Hilbert-Schmidt independence criterion \citep[HSIC,][]{gretton05a} which will also be discussed. Finally, we will review feature selection techniques introduced in machine learning which make use of these dependence measures.

\subsection{Mutual information}

Mutual information \citep[MI,][]{shannon48} is a symmetric measure of dependence which is related to the entropy.
Assuming $X$ and $Y$ are two random vectors which are absolutely continuous with respect to the Lebesgue measure on $\mathbb{R}^p$ and $\mathbb{R}^q$ with density functions $p_X(x)$ and $p_Y(y)$, respectively, one can define their marginal entropy:
\begin{equation*}
H(X) = - \int_{\mathbb{R}^p} p_X(x) \ln(p_X(x)) dx
\end{equation*}
and $H(Y)$ similarly. Denoting $p_{X,Y}(x,y)$ their joint density function, the joint entropy between $X$ and $Y$ writes
\begin{equation*}
H(X,Y) = - \int_{\mathbb{R}^{p+q}} p_{X,Y}(x,y) \ln(p_{X,Y}(x,y)) dxdy.
\end{equation*}
MI is then formally defined as
\begin{align*}
I(X;Y) &= H(X) + H(Y) - H(X,Y) \\
&=  \int_{\mathbb{R}^{p+q}} p_{X,Y}(x,y) \ln\left(\frac{p_{X,Y}(x,y)}{p_Y(y)p_{X}(x)}\right)dxdy.
\end{align*}
Interestingly, MI equals zero if and only if $X$ and $Y$ are independent. This implies that MI is able to detect nonlinear 
dependencies between random variables, unlike the correlation coefficient. 
It is also easy to check that $I(X;Y) \geq 0$ with Jensen's inequality. Further note that it is not a distance since it does not obey 
the triangle inequality. A simple modified version yielding a distance, the variation of information (VI), is given by
\begin{equation*}
VI(X;Y) = H(X,Y) - I(X;Y) = H(X) + H(Y) - 2I(X;Y).
\end{equation*}
Another variant is the squared-loss mutual information \citep[SMI,][]{suzuki09}:
\begin{equation}
SMI(X;Y) = \int_{\mathbb{R}^{p+q}} \left(\frac{p_{X,Y}(x,y)}{p_Y(y)p_{X}(x)}-1\right)^2 p_Y(y)p_{X}(x)dxdy \label{SMI}
\end{equation}
which is again a dependence measure verifying $SMI(X;Y)\geq0$ with equality if and only if $X$ and $Y$ are independent.
Applications of MI, VI and SMI include independent component analysis \citep{hyva00}, image registration \citep{pluim03} and hierarchical clustering \citep{meilua07}, among many others.\\

In the context of global sensitivity analysis, we have seen in Section \ref{fdiv} that MI and SMI arise as sensitivity indices when specific 
Csisz\'ar f-divergences are chosen to evaluate the dissimilarity between the output $Y$ and the conditional output $Y|X^k$ where $X^k$ is one of the input variables. We will then study the two following sensitivity indices:
\begin{equation}
S_{X^k}^{MI} = I(X^k;Y) = \int_{\mathbb{R}^2} p_{X^k,Y}(x,y) \ln\left(\frac{p_{X^k,Y}(x,y)}{p_Y(y)p_{X^k}(x)}\right)dxdy \label{Smi}
\end{equation}
and
\begin{equation}
S_{X^k}^{SMI} = SMI(X^k;Y) = \int_{\mathbb{R}^2} \left(\frac{p_{X^k,Y}(x,y)}{p_Y(y)p_{X^k}(x)}-1\right)^2 p_Y(y)p_{X^k}(x)dxdy. \label{Ssmi}
\end{equation}
A normalized version of $S_{X^k}^{MI}$ given by $I(X^k;Y)/H(Y)$  has already been proposed by \citet{krzy01} and compared to Sobol sensitivity indices by \citet{auder08}.

\subsection{Distance correlation}

The distance correlation was introduced by \citet{dcor07} to address the problem of testing dependence between two random vectors $X$ in $\mathbb{R}^p$ and $Y$ in $\mathbb{R}^q$. It is based on the concept of distance covariance which measures the distance between the joint 
characteristic function of $(X,Y)$ and the product of the marginal characteristic functions.

More precisely, denote $\phi_X$ and $\phi_Y$ the characteristic function of $X$ and $Y$, respectively, and $\phi_{X,Y}$ their joint 
characteristic function. For a complex-valued function $\phi(\cdot)$, we also denote $\bar{\phi}$ the complex conjugate of $\phi$ and 
$\lvert \phi\rvert^2 = \phi\bar{\phi}$. The distance covariance (dCov) $\mathcal{V}(X,Y)$ between $X$ and $Y$ with finite first moment is then defined as a weighted $L_2$-distance between $\phi_{X,Y}$ and $\phi_X\phi_Y$ given by
\begin{align}
\mathcal{V}^2(X,Y) &= \lVert \phi_{X,Y} - \phi_X\phi_Y \rVert_w^2 \nonumber \\
&= \int_{\mathbb{R}^{p+q}} \lvert \phi_{X,Y} (t,s) - \phi_X(t)\phi_Y(s)\rvert^2 w(t,s)dtds \label{dcov}
\end{align}
where the weight function $w(t,s) = (c_pc_q \lVert t\rVert_2^{1+p}\lVert s\rVert_2^{1+q})^{-1}$ with constants $c_l = \pi^{(1+l)/2}/\Gamma((1+l)/2)$ for $l\in\mathbb{N}$ is chosen to ensure invariance properties, see \citet{dcor07}. The distance correlation (dCor) $\mathcal{R}(X,Y)$ between $X$ and $Y$ is then naturally defined as
\begin{equation}
\mathcal{R}^2(X,Y) = \frac{\mathcal{V}^2(X,Y)}{\sqrt{\mathcal{V}^2(X,X)\mathcal{V}^2(Y,Y)}} \label{dcor}
\end{equation}
if $\mathcal{V}^2(X,X)\mathcal{V}^2(Y,Y) > 0$ and equals $0$ otherwise. Important properties of the distance correlation introduced in (\ref{dcor}) include that $0\leq\mathcal{R}(X,Y)\leq 1$ and $\mathcal{R}(X,Y)=0$ if and only if $X$ and $Y$ are independent. Interestingly, the distance covariance in (\ref{dcov}) can be computed in terms of expectations of pairwise Euclidean distances, namely
\begin{align}
\mathcal{V}^2(X,Y) &= \mathbb{E}_{X,X',Y,Y'} \lVert X-X' \rVert_2 \lVert Y-Y' \rVert_2 \nonumber \\
&+ \mathbb{E}_{X,X'} \lVert X-X' \rVert_2 \mathbb{E}_{Y,Y'} \lVert Y-Y' \rVert_2 \nonumber \\
&-2 \mathbb{E}_{X,Y} \left[ \mathbb{E}_{X'} \lVert X-X' \rVert_2 \mathbb{E}_{Y'} \lVert Y-Y' \rVert_2\right] \label{dcovE}
\end{align}
where $(X',Y')$ is an i.i.d. copy of $(X,Y)$. Concerning estimation, let $(X_i,Y_i)_{i=1,\ldots,n}$ be a sample of the random vector $(X,Y)$. Following equation (\ref{dcovE}), an estimator $\mathcal{V}^2_n(X,Y)$ of $\mathcal{V}^2(X,Y)$ is then given by
\begin{align}
\mathcal{V}^2_n(X,Y) & = \frac{1}{n^2} \sum_{i,j=1}^n \lVert X_i-X_j \rVert_2 \lVert Y_i-Y_j \rVert_2 \nonumber \\
& + \frac{1}{n^2} \sum_{i,j=1}^n \lVert X_i-X_j \rVert_2 \frac{1}{n^2} \sum_{i,j=1}^n\lVert Y_i-Y_j \rVert_2 \nonumber \\
& - \frac{2}{n} \sum_{i=1}^n \left[ \frac{1}{n} \sum_{j=1}^n\lVert X_i-X_j \rVert_2 \frac{1}{n} \sum_{j=1}^n\lVert Y_i-Y_j \rVert_2\right] . \label{dcovemp}
\end{align}
Denoting $a_{ij} = \lVert X_i-X_j \rVert_2$, $\bar{a}_{i\cdot} = \sum_j a_{ij}/n$, $\bar{a}_{\cdot j} = \sum_i a_{ij}/n$, $\bar{a}_{\cdot\cdot}=\sum_{ij}a_{ij}/n^2$, $A_{ij} = a_{ij} - \bar{a}_{i\cdot} - \bar{a}_{\cdot j}  + \bar{a}_{\cdot\cdot}$ and similarly $B_{ij} = b_{ij} - \bar{b}_{i\cdot} - \bar{b}_{\cdot j}  + \bar{b}_{\cdot\cdot}$ for $b_{ij} = \lVert Y_i-Y_j \rVert_2$, \citet{dcor07} show that equation (\ref{dcovemp}) can be written as
\begin{equation*}
\mathcal{V}^2_n(X,Y)  = \frac{1}{n^2} \sum_{i,j=1}^n A_{ij}B_{ij}
\end{equation*}
and is also equal to equation (\ref{dcov}) if one uses the empirical characteristic functions computed on the sample $(X_i,Y_i)_{i=1,\ldots,n}$. The empirical distance correlation $\mathcal{R}_n(X,Y)$ is then 
\begin{equation*}
\mathcal{R}_n^2(X,Y) = \frac{\mathcal{V}_n^2(X,Y)}{\sqrt{\mathcal{V}_n^2(X,X)\mathcal{V}^2_n(Y,Y)}}
\end{equation*}
and satisfies $0\leq\mathcal{R}_n(X,Y)\leq 1$. Although $\mathcal{V}^2_n(X,Y)$ is a consistent estimator of $\mathcal{V}^2(X,Y)$, it is easy to see that it is biased. \citet{dcor13b} propose an unbiased version of $\mathcal{V}^2_n(X,Y)$ and a specific correction for the high-dimensional case $p,q \gg 1$ is investigated in \citet{dcor13a}.
Further note that \citet{dcor07} also study $\mathcal{V}^{2(\alpha)}(X,Y)$ defined as 
\begin{align}
\mathcal{V}^{2(\alpha)}(X,Y) &= \int_{\mathbb{R}^{p+q}} \lvert \phi_{X,Y} (t,s) - \phi_X(t)\phi_Y(s)\rvert^2 w_{\alpha}(t,s)dtds \nonumber \\
&= \mathbb{E}_{X,X',Y,Y'} \lVert X-X' \rVert_2^{\alpha} \lVert Y-Y' \rVert_2^{\alpha} \nonumber \\
&+ \mathbb{E}_{X,X'} \lVert X-X' \rVert_2^{\alpha} \mathbb{E}_{Y,Y'} \lVert Y-Y' \rVert_2^{\alpha} \nonumber \\
&-2 \mathbb{E}_{X,Y} \left[ \mathbb{E}_{X'} \lVert X-X' \rVert_2^{\alpha} \mathbb{E}_{Y'} \lVert Y-Y' \rVert_2^{\alpha}\right] \label{dcova}
\end{align}
with the new weight function $w_{\alpha}(t,s) = (C(p,\alpha)C(q,\alpha) \lVert t\rVert_2^{\alpha+p}\lVert s\rVert_2^{\alpha+q})^{-1}$ and constants $C(l,\alpha) = \frac{2\pi^{l/2}\Gamma(1-\alpha/2)}{\alpha 2^{\alpha}\Gamma((l+\alpha)/2)}$ as soon as $\mathbb{E}(\lVert X\rVert_2^{\alpha} + \lVert Y\rVert_2^{\alpha} ) < \infty$ and $0 < \alpha < 2$. Distance covariance is retrieved for $\alpha=1$. The very general case of $X$ and $Y$ living in metric spaces has been examined by \citet{lyons13}. More precisely, let $(\mathcal{X},\rho_\mathcal{X})$ and $(\mathcal{Y},\rho_\mathcal{Y})$ be metric spaces of negative type \citep[see][]{lyons13}, the generalized distance covariance 
\begin{align}
\mathcal{V}^{2}_{\rho_\mathcal{X},\rho_\mathcal{Y}}(X,Y) & = \mathbb{E}_{X,X',Y,Y'} \rho_\mathcal{X}(X,X')\rho_\mathcal{Y}(Y,Y' )\nonumber \\
&+ \mathbb{E}_{X,X'} \rho_\mathcal{X}(X,X') \mathbb{E}_{Y,Y'} \rho_\mathcal{Y}(Y,Y' ) \nonumber \\
&-2 \mathbb{E}_{X,Y} \left[ \mathbb{E}_{X'} \rho_\mathcal{X}(X,X') \mathbb{E}_{Y'} \rho_\mathcal{Y}(Y,Y' )\right] \label{dcovmetric}
\end{align}
characterizes independence between $X\in\mathcal{X}$ and $Y\in\mathcal{Y}$.\\

Coming back to sensitivity analysis, just like we defined a new index based on mutual information, we can finally introduce an index based on distance correlation, i.e.
\begin{equation}
S_{X^k}^{dCor} = \mathcal{R}(X^k,Y) \label{Sdcor}
\end{equation}
which will measure the dependence between an input variable $X^k$ and the output $Y$. Since distance correlation is designed to detect nonlinear relationships, we except this index to quantify effectively the impact of $X^k$ on $Y$. Besides, considering that distance covariance is defined in arbitrary dimension, this index generalizes easily to the multivariate case:
\begin{equation*}
S_{X^u}^{dCor} = \mathcal{R}(X^u,Y)
\end{equation*}
for evaluating the impact of a group of inputs $X^u$ on a multivariate output $Y$.

\begin{remark}
\label{remPF}
The limiting case $\alpha \rightarrow 2$ in (\ref{dcova}) interestingly leads to $\mathcal{V}^{2(2)}(X,Y) = \textrm{Cov}(X,Y)^2$, see \citet{dcor07}.
This turns out to be another original way for defining a new sensitivity index. Indeed, recall that Sobol first-oder sensitivity index actually equals $\textrm{Cov}(Y,Y_{X^k})/\textrm{Var}(Y)$ where $Y_{X^k}$ is an independent copy of $Y$ obtained by fixing $X^k$, see \citet{janon13}. The idea is then to replace the covariance (obtained with $\alpha \rightarrow 2$) by dCov (with $\alpha=1$):
\begin{equation}
S_{X^k}^{dCorPF} = \mathcal{R}(Y,Y_{X^k}),  \label{SdcorPF}
\end{equation}
where PF stands for pick-and-freeze, since this index generalizes the pick-and-freeze estimator proposed by \citet{janon13} and is able to detect nonlinear dependencies, unlike the correlation coefficient.
\end{remark}

\subsection{HSIC} \label{secHSIC}

\subsubsection{Definition}

The Hilbert-Schmidt independence criterion proposed by \citet{gretton05a} builds upon kernel-based approaches for detecting dependence, and more particularly on cross-covariance operators in RKHSs. Here, we only give a brief summary and introduction on this topic and refer the reader to \citet{berlinet04,gretton05a,smola07} for details.

Let the random vector $X\in\mathcal{X}$ have distribution $P_X$ and consider a RKHS $\mathcal{F}$ of functions $\mathcal{X}\rightarrow \mathbb{R}$ with kernel $k_\mathcal{X}$ and dot product $\langle\cdot,\cdot\rangle_\mathcal{F}$. Similarly, we can also define a second RKHS  $\mathcal{G}$ of functions $\mathcal{Y}\rightarrow \mathbb{R}$ with kernel $k_\mathcal{Y}$ and dot product $\langle\cdot,\cdot\rangle_\mathcal{G}$ associated to the random vector $Y\in\mathcal{Y}$ with distribution $P_Y$. By definition, the cross-covariance operator $C_{XY}$ associated to the joint distribution $P_{XY}$ of $(X,Y)$ is the linear operator $\mathcal{G}\rightarrow \mathcal{F}$ defined for every $f\in\mathcal{F}$ and $g\in\mathcal{G}$ as
\begin{equation*}
\langle f,C_{XY}g\rangle_\mathcal{F} = \mathbb{E}_{XY}[f(X)g(Y)] -\mathbb{E}_Xf(X)\mathbb{E}_Yg(Y).
\end{equation*}
In a nutshell, the cross-covariance operator generalizes the covariance matrix by representing higher order correlations between $X$ and $Y$ through nonlinear kernels.
For every linear operator $C:\mathcal{G}\rightarrow \mathcal{F}$ and provided the sum converges, the Hilbert-Schmidt norm of $C$ is given by
\begin{equation*}
\lVert C\rVert_{HS}^2 = \sum_{k,l} \langle u_k,Cv_l\rangle_\mathcal{F} 
\end{equation*}
where $u_k$ and $v_l$ are orthonormal bases of $\mathcal{F}$ and $\mathcal{G}$, respectively. This is simply the generalization of the Frobenius norm on matrices. The HSIC criterion is then defined as the Hilbert-Schmidt norm of the cross-covariance operator:
\begin{align}
HSIC(X,Y)_{\mathcal{F},\mathcal{G}} & = \lVert C_{XY}\rVert_{HS}^2 \nonumber \\
& = \mathbb{E}_{X,X',Y,Y'} k_\mathcal{X}(X,X')k_\mathcal{Y}(Y,Y' )\nonumber \\
&+ \mathbb{E}_{X,X'} k_\mathcal{X}(X,X') \mathbb{E}_{Y,Y'} k_\mathcal{Y}(Y,Y' ) \nonumber \\
&-2 \mathbb{E}_{X,Y} \left[ \mathbb{E}_{X'} k_\mathcal{X}(X,X') \mathbb{E}_{Y'} k_\mathcal{Y}(Y,Y' )\right] \label{HSIC}
\end{align}
where the last equality in terms of kernels is proven in \citet{gretton05a}. An important property of $HSIC(X,Y)_{\mathcal{F},\mathcal{G}}$ is that it equals $0$ if and only if $X$ and $Y$ are independent, as long as the associated RKHSs $\mathcal{F}$ and $\mathcal{G}$ are universal, i.e. they are dense in the space of continuous functions with respect to the infinity norm \citep{gretton05b}. Examples of kernels generating universal RKHSs are e.g. the Gaussian and the Laplace kernels \citep{sri09}.

It is interesting to note the similarity between the generalized distance covariance of equation (\ref{dcovmetric}) and the HSIC criterion (\ref{HSIC}). Actually, \citet{sejdi13} recently studied the deep connection between these approaches and show that 
\begin{equation*}
\mathcal{V}^{2}_{\rho_\mathcal{X},\rho_\mathcal{Y}}(X,Y) = 4\ HSIC(X,Y)_{\mathcal{F},\mathcal{G}}
\end{equation*}
if the kernels $k_\mathcal{X}$ and $k_\mathcal{Y}$ generate the metrics $\rho_\mathcal{X}$ and $\rho_\mathcal{Y}$, respectively \citep[see][]{sejdi13}. In particular, the standard distance covariance (\ref{dcovE}) is retrieved with the (universal) kernel $k(z,z') = \frac{1}{2}(\lVert z \rVert_2 + \lVert z' \rVert_2 - \lVert z - z'\rVert_2 )$ which generates the metric $\rho(z,z') = \lVert z - z'\rVert_2$.

Assume now that $(X_i,Y_i)_{i=1,\ldots,n}$ is a sample of the random vector $(X,Y)$ and denote $K_\mathcal{X}$ and $K_\mathcal{Y}$ the Gram matrices with entries $K_\mathcal{X}(i,j) = k_\mathcal{X}(X_i,X_j)$ and $K_\mathcal{Y}(i,j) = k_\mathcal{Y}(Y_i,Y_j)$. \citet{gretton05a} propose the following consistent estimator for $HSIC(X,Y)_{\mathcal{F},\mathcal{G}}$:
\begin{equation*}
HSIC_n(X,Y)_{\mathcal{F},\mathcal{G}} = \frac{1}{n^2} \textrm{Tr}\left(K_\mathcal{X} H K_\mathcal{Y} H \right)
\end{equation*}
where $H$ is the centering matrix such that $H(i,j) = \delta_{ij} - \frac{1}{n}$. Besides, it is easy to check that $HSIC_n(X,Y)_{\mathcal{F},\mathcal{G}}$ can be expressed just like the empirical distance covariance (\ref{dcovemp}):
\begin{align*}
HSIC_n(X,Y)_{\mathcal{F},\mathcal{G}} & = \frac{1}{n^2} \sum_{i,j=1}^n k_\mathcal{X}(X_i,X_j) k_\mathcal{Y}(Y_i,Y_j)\\
& + \frac{1}{n^2} \sum_{i,j=1}^n  k_\mathcal{X}(X_i,X_j) \frac{1}{n^2} \sum_{i,j=1}^n k_\mathcal{Y}(Y_i,Y_j)\\
& - \frac{2}{n} \sum_{i=1}^n \left[ \frac{1}{n} \sum_{j=1}^n k_\mathcal{X}(X_i,X_j) \frac{1}{n} \sum_{j=1}^n\ k_\mathcal{Y}(Y_i,Y_j)\right].
\end{align*}
An unbiased estimator is also introduced by \citet{song12}.\\

We can finally propose a sensitivity index generalizing (\ref{Sdcor}):
\begin{equation}
S_{X^k}^{HSIC_{\mathcal{F},\mathcal{G}}} = \mathcal{R}(X^k,Y)_ {\mathcal{F},\mathcal{G}} \label{Shsic}
\end{equation}
where the kernel-based distance correlation is given by
\begin{equation*}
\mathcal{R}^2(X,Y)_ {\mathcal{F},\mathcal{G}} = \frac{HSIC(X,Y)_{\mathcal{F},\mathcal{G}} }{\sqrt{HSIC(X,X)_{\mathcal{F},\mathcal{F}} HSIC(Y,Y)_{\mathcal{G},\mathcal{G}} }}
\end{equation*}
and the kernels inducing $\mathcal{F}$ and $\mathcal{G}$ have to be chosen within the class of universal kernels. The multivariate extension of $S_{X^k}^{HSIC_{\mathcal{F},\mathcal{G}}}$ is straightforward. The impact of the choice of kernels has previously been studied by \citet{sri09} in the context of independence hypothesis tests.

\begin{remark}
Instead of working with the cross-covariance operator $C_{XY}$, \citet{fuku07} work with the normalized cross-covariance operator (NOCCO) $V_{XY}$ defined as $C_{XY} = C_{YY}^{1/2} V_{XY} C_{XX}^{1/2}$, see \citet{fuku07} for the existence of this representation. Just as the HSIC criterion, the associated measure of dependence is given by $I^{NOCCO}(X,Y)=\lVert V_{XY} \rVert_{HS}^2$. Interestingly, $I^{NOCCO}(X,Y)$ is independent of the choice of kernels and is actually equal to the squared-loss mutual information (\ref{SMI}) under some assumptions, see \citet{fuku08}. Despite the advantage of being kernel-free, using $I^{NOCCO}$ in practice unfortunately requires to work with an estimator with a regularization parameter, which has to be selected \citep{fuku07}. Nevertheless, it is still interesting to use this approach for approximating SMI efficiently, since dimensionality limitations related to density function estimation no longer apply.
\end{remark}

\begin{remark}
The pick-and-freeze estimator defined in Remark \ref{remPF} can be readily generalized with kernels:
\begin{equation}
S_{X^k}^{HSIC_{\mathcal{F},\mathcal{G}}PF} = \mathcal{R}(Y,Y_{X^k})_ {\mathcal{G},\mathcal{G}} \label{ShsicPF}
\end{equation}
where this time only the kernel acting on $\mathcal{Y}$ needs to be specified.
\end{remark}

\subsubsection{Going beyond $Y\in\mathbb{R}^q$} \label{secFunc}

The kernel point of view in HSIC also provides an elegant and powerful framework for dealing with categorical inputs and outputs, as well as functional ones.
 
The categorical case is common practice in feature selection, since the target output is often represented as labels. Appropriate kernels include for example $k_\mathcal{Y}(y,y') = \delta_{yy'}/n_y$ where $n_y$ is the number of samples with label $y$, see e.g. \citet{song12,yamada13}.
From a GSA perspective, this implies that we can evaluate the impact of the inputs on level sets of the output by a simple change of variable 
$Z =  \mathds{1}\{Y > t\}$ for a given threshold $t$. We can note the resemblance with the approach of \citet{fort13} if one uses a contrast function adapted to exceedance probabilities.

As a matter of fact, it is also possible to design dedicated semi-metrics for functional data which can be incorporated in the definition of the kernels, see e.g. \citet{ferraty06}.
For example, let $\Delta(\cdot,\cdot)$ be such a semi-metric defined on $\mathcal{Y}\times\mathcal{Y}$ when the output variable is of functional type. The kernel associated to $\mathcal{Y}$ is then given by $k_\mathcal{Y}(y,y') = k(\Delta(y,y'))$ where $k$ is a kernel acting on $\mathbb{R}$. The same scheme applies to functional inputs as well, see \citet{gin12} for an illustration in the context of surrogate modeling where the semi-metric is a cheap and simplified computer code. However, a theoretical shortcoming lies in our current inability to check if such semi-metric kernels are universal, which implies that we can not claim that independence can be detected. Despite this deficiency, we show in Section \ref{secexp} that from a practical perspective, the use of a semi-metric based on principal components can efficiently deal with a functional output given as a 2D map.

\subsection{Feature selection as an alternative to screening} \label{secfs}

In machine learning, feature selection aims at identifying relevant features (among a large set) with respect to a prediction task. The goal is to detect irrelevant or redundant features which may increase the prediction variance without reducing its bias. As a matter of fact, this closely resembles the objective of factor screening in GSA. The main difference is that in GSA, input variables are usually assumed to be independent, whereas in feature selection redundant features, i.e. highly dependent factors, precisely have to be filtered out. This apparently naive distinction actually makes feature selection an interesting alternative to screening when some input variables are correlated. But it is important to note that it is also a powerful option even in the independent case.
We do not intend here to give an exhaustive review of feature selection techniques, but rather detail some approaches which make use of the dependence measures we recapped above. We hope that it will illustrate how they can be used as new screening procedures in high dimensional problems.\\

Literature on feature selection is abundant and entails many approaches. In the high dimensional setting, model-based techniques include for example the Lasso \citep{tib96} or sparse additive models \citep{ravi09}, see \citet{fan10} for a selective overview. Generalizations for the ultra-high dimensional case usually replace penalty-based techniques to focus on marginal regression, where an underlying model is still assumed (e.g. linear \citet{fan08} or non-parametric additive \citet{fan11}). Another line of work for the ultra-high dimensional setting are model-free methods, where only dependence measures are used to identify relevant features. Except for the very specific HSIC Lasso technique \citep{yamada13}, here we only focus on pure dependence-based approaches. 

Let us first introduce the concept of Max-Dependency \citep{peng05}. Denote $X^1,\ldots,X^p$ the set of available features, $Y$ the target output to predict and $D(\cdot,\cdot)$ any measure quantifying the dependence between two random vectors. The Max-Dependency scheme for feature selection involves finding $m$ features $X^{i_1},\ldots,X^{i_m}$ which jointly have the largest dependency with $Y$, i.e. one has to solve the following optimization problem
\begin{equation}
\max_{ \left\{i_1,\ldots,i_m\right\} \subset \left\{1,\ldots,p\right\} } D(\left\{X^{i_1},\ldots,X^{i_m}\right\},Y). \label{MaxD}
\end{equation}
Solving (\ref{MaxD}) is however computationally infeasible when $m$ and $p$ are large for cardinality reasons. Near-optimal solutions are then usually found by iterative procedures, where features are added one at a time in the subset $X^{i_1},\ldots,X^{i_m}$ (forward selection). On the other hand, the dependence measure $D(\cdot,\cdot)$ must also be robust to dimensionality, which is hard to achieve in practice when the number of samples is less than $m$. Consequently, marginal computations which only involve $D(X^k,Y)$ terms are usually preferred. The Max-Relevance criterion \citep{peng05} serves in this context as a proxy to Max-Dependency, where the optimization problem writes
\begin{equation}
\max_{ \left\{i_1,\ldots,i_m\right\} \subset \left\{1,\ldots,p\right\} } \frac{1}{m} \sum_{k=1}^m D(X^{i_k},Y). \label{MaxR}
\end{equation}
But when the features are dependent, it is likely that this criterion will select redundant features. To limit this effect, one can add a condition of Min-Redudancy expressed as
\begin{equation}
\min_{ \left\{i_1,\ldots,i_m\right\} \subset \left\{1,\ldots,p\right\} } \frac{1}{m^2} \sum_{k,l=1}^m D(X^{i_k},X^{i_l}). \label{MinR}
\end{equation}
The final scheme combining (\ref{MaxR}) and (\ref{MinR}), called minimal-redundancy-maximal-relevance (mRMR), is given by
\begin{equation}
\max_{ \left\{i_1,\ldots,i_m\right\} \subset \left\{1,\ldots,p\right\} } \frac{1}{m} \sum_{k=1}^m D(X^{i_k},Y) -  \frac{1}{m^2} \sum_{k,l=1}^m D(X^{i_k},X^{i_l}) . \label{mrmr}
\end{equation}
Forward and backward procedures for mRMR are investigated by \citet{peng05} where $D(\cdot,\cdot)$ is chosen as the mutual information. Similarly, forward and backward approaches where MI is replaced with HSIC is introduced by \citet{song12}.
A purely marginal point of view is studied by \citet{li12} where the authors propose the dCor criterion (\ref{dcor}). In a nutshell, the dCor measure is computed between $Y$ and each factor $X^k$, $k=1,\ldots,p$ and only the features with dCor above a certain threshold are retained. A sure screening property of this approach is also proven. \citet{bala13} extend this work by considering a modified version of the HSIC dependence measure (supremum of HSIC over a family of universal kernels, denoted sup-HSIC). Even if the sure screening procedure of this generalized method is proven, the authors mention that every feature selection technique based on marginal computations fails at detecting features that may be marginally uncorrelated with the output but are in fact jointly correlated with it. As a result, they propose the following iterative approach:
\begin{enumerate}
\item Compute the marginal sup-HSIC measures between $Y$ and each feature $X^k$, $k=1,\ldots,p$ and select the inputs with a measure above a given threshold. Let $u\subset \left\{1,\ldots,p\right\} $ be the subset of selected features.
\item Compute sup-HSIC between $Y$ and $\left(X^{u},X^k\right)$ for each $k \notin u$. Augment $u$ with features having a measure greater than the sup-HSIC criterion between $Y$ and $X^{u}$.
\item Repeat until the subset of selected features stabilizes or when its cardinality reaches a given maximum value.
\end{enumerate}
As pointed out previously, another drawback of marginal computations which is not taken care of by the above scheme is that redundant variables are not eliminated. But \citet{bala13} design another iterative procedure to deal with this case.
Finally, let us note that in the examples of Section \ref{secexp}, we will only study the above iterative technique since we focus on independent input factors. We plan to investigate in particular the full mRMR approach for problems with correlated inputs in a future work.

Instead of working with forward and backward approaches, \citet{yamada13} propose a combination of the Lasso and the HSIC dependence measure. Denote $\tilde{K}^k_\mathcal{X} = HK^k_\mathcal{X}H$ for $k=1,\ldots,p$ and $\tilde{K}_\mathcal{Y} = HK_\mathcal{Y}H$ the centered Gram matrices computed from a sample $\left(X^1_i,\ldots,X^p_i,Y_i\right)_{i=1,\ldots,n}$ of $\left(X^1,\ldots,X^p,Y\right)$ following the notations of Section \ref{secHSIC}. The HSIC Lasso solves the following optimization problem
\begin{equation}
\min_{\alpha\in \mathbb{R}^p} \frac{1}{2} \lVert \tilde{K}_\mathcal{Y} -\sum_{k=1}^p \alpha_k \tilde{K}^k_\mathcal{X}\rVert^2_{\textrm{Frob}} + \lambda \lVert \alpha \rVert_1 \label{hsiclasso}
\end{equation}
with constraints $\alpha_1,\ldots,\alpha_p \geq 0$ and where $\lVert \cdot \rVert_{\textrm{Frob}}$ stands for the Frobenius norm and $\lambda$ is a regularization parameter. Interestingly, the first term of equation (\ref{hsiclasso}) expands as
\begin{align*}
\frac{1}{2n^2} \lVert \tilde{K}_\mathcal{Y} -\sum_{k=1}^p \alpha_k \tilde{K}^k_\mathcal{X}\rVert^2_{\textrm{Frob}} &= 
\frac{1}{2} HSIC_n(Y,Y) - \sum_{k=1}^p \alpha_k HSIC_n(X^k,Y) \\
& + \frac{1}{2} \sum_{k,l=1}^p \alpha_k \alpha_l HSIC_n(X^k,X^l)
\end{align*}
using that $\tilde{K}^k_\mathcal{X}$, $\tilde{K}_\mathcal{Y}$ are symmetric and $H$ is idempotent, which highlights the strong correspondence with the mRMR criterion (\ref{mrmr}). The authors show that (\ref{hsiclasso}) can be recast as a standard Lasso program and propose a dual augmented Lagrangian algorithm to solve the optimization problem. They also discuss a variant based on the $I^{NOCCO}$ dependence measure.

\begin{remark}
We mentioned before that feature selection techniques based on dependence measures have been particularly designed for the ultra-high dimensional case, which is not the common setting of screening problems in GSA. Nevertheless, we illustrate in Section \ref{secexp} that they perform remarkably well on complex benchmark functions, while requiring very few samples of the output variable. This reveals their high potential for preliminary screening of expensive computer codes.
\end{remark}

\section{Experiments} \label{secexp}

In this Section, we finally assess the performance of all the new sensitivity indices introduced before on a series of benchmark analytical functions and two industrial applications. All benchmark functions can be found in the Virtual Library of Simulation Experiments available at \url{http://www.sfu.ca/~ssurjano/index.html}. For easier comparison, we first summarize the proposed indices in Table \ref{tabSI} (SI stands for sensitivity index and see \citet{tarantola06} for RBD-FAST).
\begin{table}[h]
   \centering
   \scalebox{0.8}{
   \begin{tabular}{| c | c | c |}
   \hline
         & & \\
      Index & Origin & Notes \\
            & & \\
      \hline
      & & \\
      $S_{X}^1$ & Sobol first-order SI & Normalized version (RBD-FAST est.)\\
      & & \\
      \hline
      & & \\
      $S_{X}^{TOT}$ & Sobol total SI & Normalized version (RBD-FAST est.)\\
      & & \\
      \hline
      & & \\
      $S_{X}^f$ (eq. (\ref{Sif})) & Csisz\'ar f-divergences & Includes as special case \citet{borgonovo07}, \\
      & & $S_{X}^{MI}$ (eq. (\ref{Smi})) and $S_{X}^{SMI}$ (eq. (\ref{Ssmi})) \\
      & & \\
      \hline
      & & \\
      $S_{X}^{dCor}$ (eq. (\ref{Sdcor})) & Distance correlation & \\
      & & \\
      \hline
      & & \\
      $S_{X}^{dCorPF}$ (eq. (\ref{SdcorPF})) & Pick-and-freeze distance correlation &  Generalization of \citet{janon13} \\
      & & \\
      \hline
      & & \\
      $S_{X}^{HSIC_{\mathcal{F},\mathcal{G}}}$ (eq. (\ref{Shsic}) & HSIC & Can be extended to categorical/functional data \\
      & & \\
      \hline
      & & \\
      $S_{X}^{HSIC_{\mathcal{F},\mathcal{G}}PF}$  (eq. (\ref{ShsicPF}) & Pick-and-freeze HSIC & Generalization of \citet{janon13} \\
      & & \\
      \hline
   \end{tabular}}
   \caption{Summary of sensitivity indices.} \label{tabSI}
\end{table}

\subsection{Analytical functions}

\paragraph{Standard GSA.} For the first experiments, we focus on GSA problems where the dimensionality is not too large (less than $10$ input variables). The objective is to compare the information given by the new indices with Sobol first-order and total indices.
\begin{itemize}
\item \citet{link06} decreasing function 
\[\eta_1(X) = 0.2X_1 + \frac{0.2}{2}X_2 + \frac{0.2}{4}X_3 + \frac{0.2}{8}X_4 + \frac{0.2}{16}X_5 + \frac{0.2}{32}X_6 + \frac{0.2}{64}X_7 + \frac{0.2}{128}X_8\]
with $X_i\sim \mathcal{U}(0,1)$, $i=1,\ldots,10$.
\end{itemize}

\begin{figure}[h]
\centering
\includegraphics[width=0.9\textwidth]{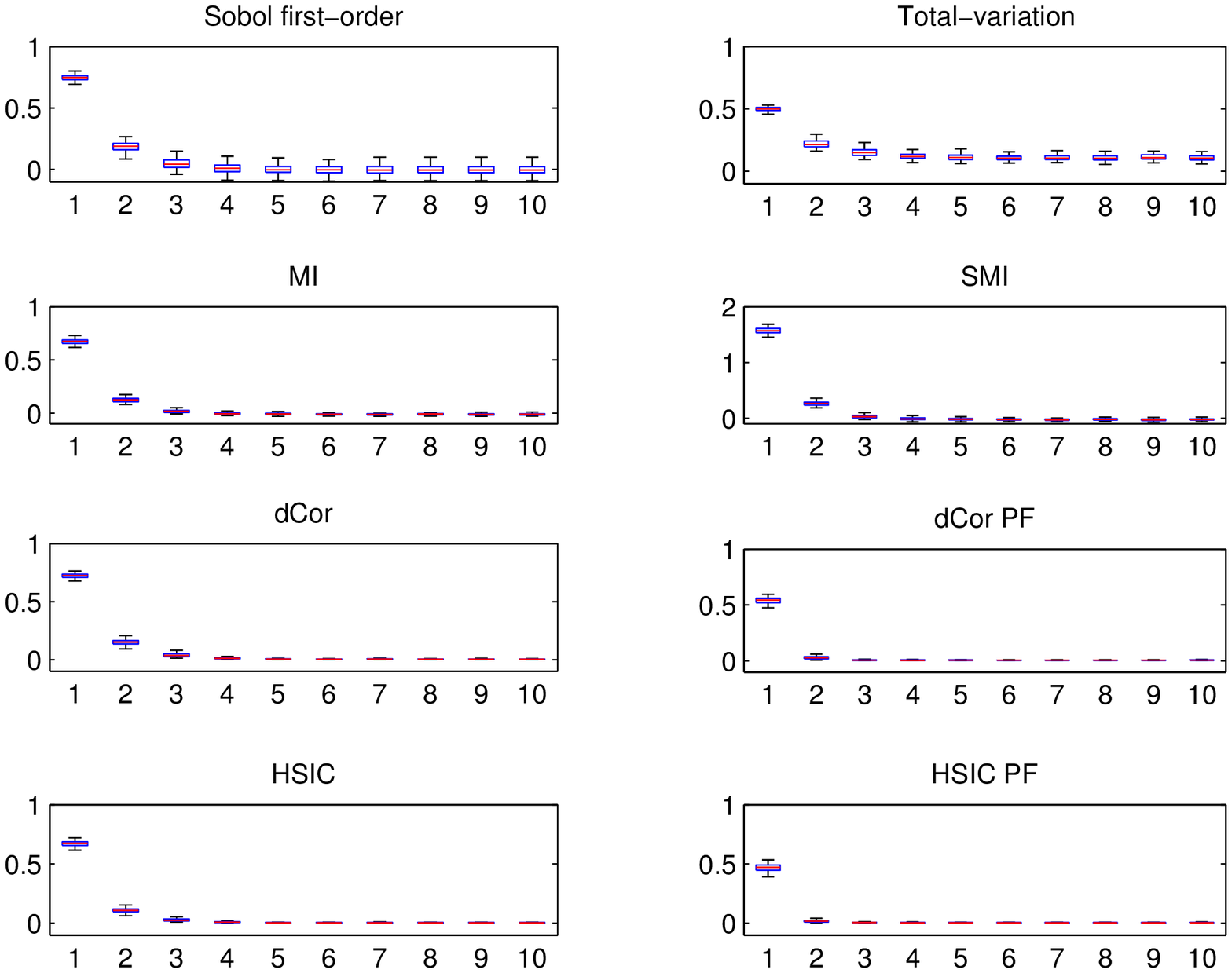}
\caption{First-order SI, $S_{X}^f$, $S_{X}^{dCor}$, $S_{X}^{dCorPF}$, $S_{X}^{HSIC_{\mathcal{F},\mathcal{G}}}$ and $S_{X}^{HSIC_{\mathcal{F},\mathcal{G}}PF}$ for function $\eta_1$, $n=500$, 100 replicates.} \label{figLink}
\end{figure}

We compute the sensitivity indices based on Csisz\'ar f-divergences, dCor, pick-and-freeze dCor, HSIC and pick-and-freeze HSIC (Gaussian kernels) with a sample of size $n=500$  and we repeat this calculation 100 times. Here we use a simple kernel density estimator since we only study first-order indices. Results are given in Figure \ref{figLink}.
Analytical first-order SIs are $S_{X_k}^1 = S_{X_k}^{TOT} = \frac{3}{4} \left(\frac{1}{4}\right)^{i-1}/(1-(\frac{1}{4})^{10})$ for $i=1,\ldots,10$, which is coherent with the estimates at the top left. As expected, indices given by Csisz\'ar f-divergences, $S_{X}^{dCor}$, $S_{X}^{dCorPF}$, $S_{X}^{HSIC_{\mathcal{F},\mathcal{G}}}$ and $S_{X}^{HSIC_{\mathcal{F},\mathcal{G}}PF}$ provide the same information as variance-based ones in this simple case of a linear model. However, let us note that dCor and HSIC detect non-influential factors very easily and robustly.

\begin{itemize}
\item \citet{loep13} function 
\[\eta_2(X) = 6X_1 + 4X_2+5.5X_3+3X_1X_2+2.2X_1X_3+1.4X_2X_3+X_4+0.5X_5+0.2X_6+0.1X_7\]
 with  $X_i\sim \mathcal{U}(0,1)$, $i=1,\ldots,10$ (the original function has constraint $\sum_{i=1}^{10} X_i=1$ but we do not consider it here).
\end{itemize}

\begin{figure}[h]
\centering
\includegraphics[width=0.9\textwidth]{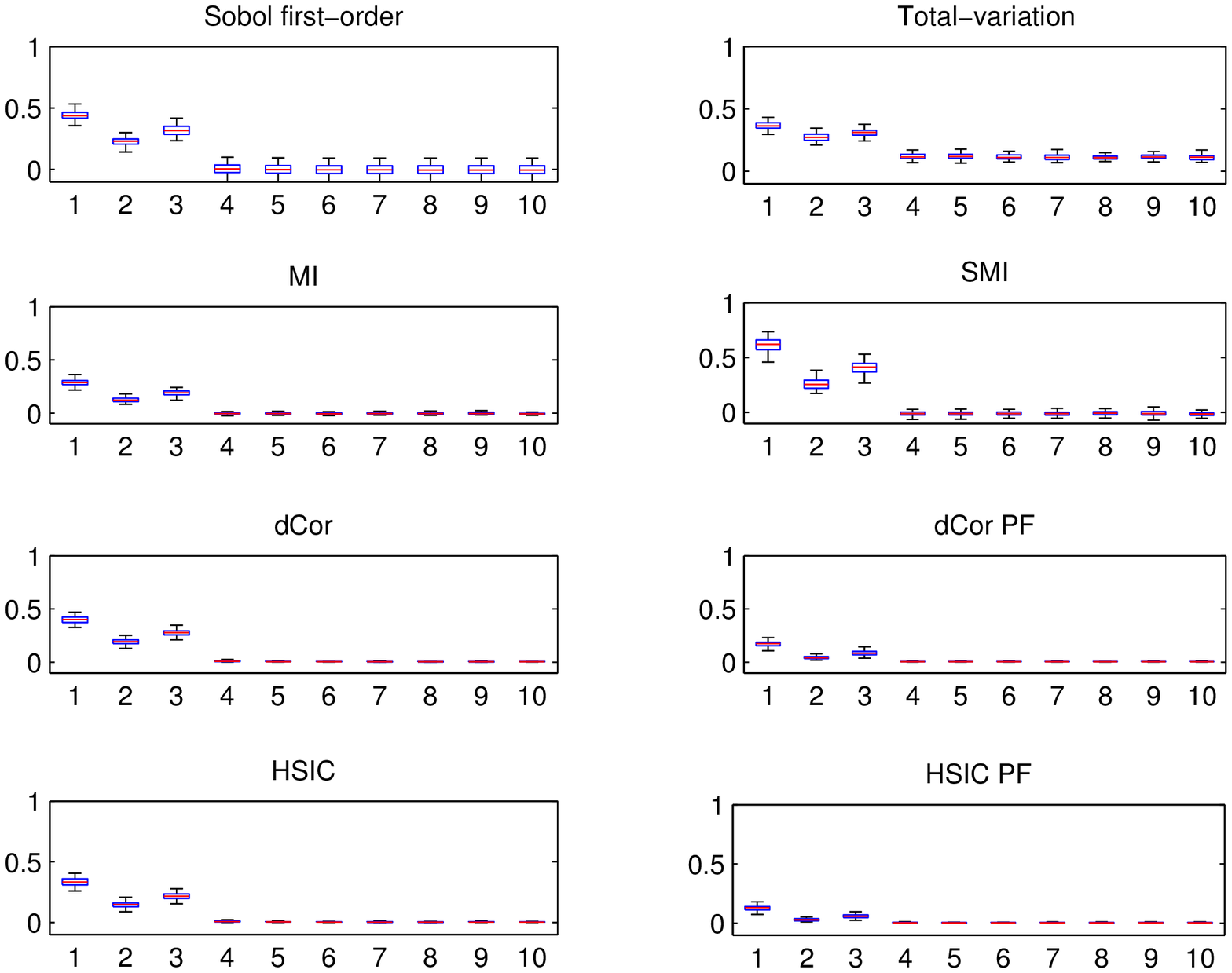}
\caption{First-order SI, $S_{X}^f$, $S_{X}^{dCor}$, $S_{X}^{dCorPF}$, $S_{X}^{HSIC_{\mathcal{F},\mathcal{G}}}$ and $S_{X}^{HSIC_{\mathcal{F},\mathcal{G}}PF}$ for function $\eta_2$, $n=500$, 100 replicates.} \label{figLoep}
\end{figure}

Conclusions are similar for the \citet{loep13} function, where only the first three inputs have a large impact on the output with very small interactions (total SIs almost equal first-order SIs), see Figure \ref{figLoep}. Note that $S_{X}^{dCorPF}$ and $S_{X}^{HSIC_{\mathcal{F},\mathcal{G}}PF}$ recover $S_X^1$ since interactions are small. Again, dCor and HSIC clearly identify inputs which are independent of the output.

\begin{itemize}
\item Ishigami function  \citep{ishigami90}
\[\eta_3(X) = \sin(X_1) + 5\sin^2(X_2) + 0.1 (X_3)^4 \sin(X_1)\]
with $X_i\sim \mathcal{U}(-\pi,\pi)$, $i=1,\ldots,3$ and constants taken from \citet{borgonovo07}.
\end{itemize}

\begin{figure}[h]
\centering
\includegraphics[width=0.9\textwidth]{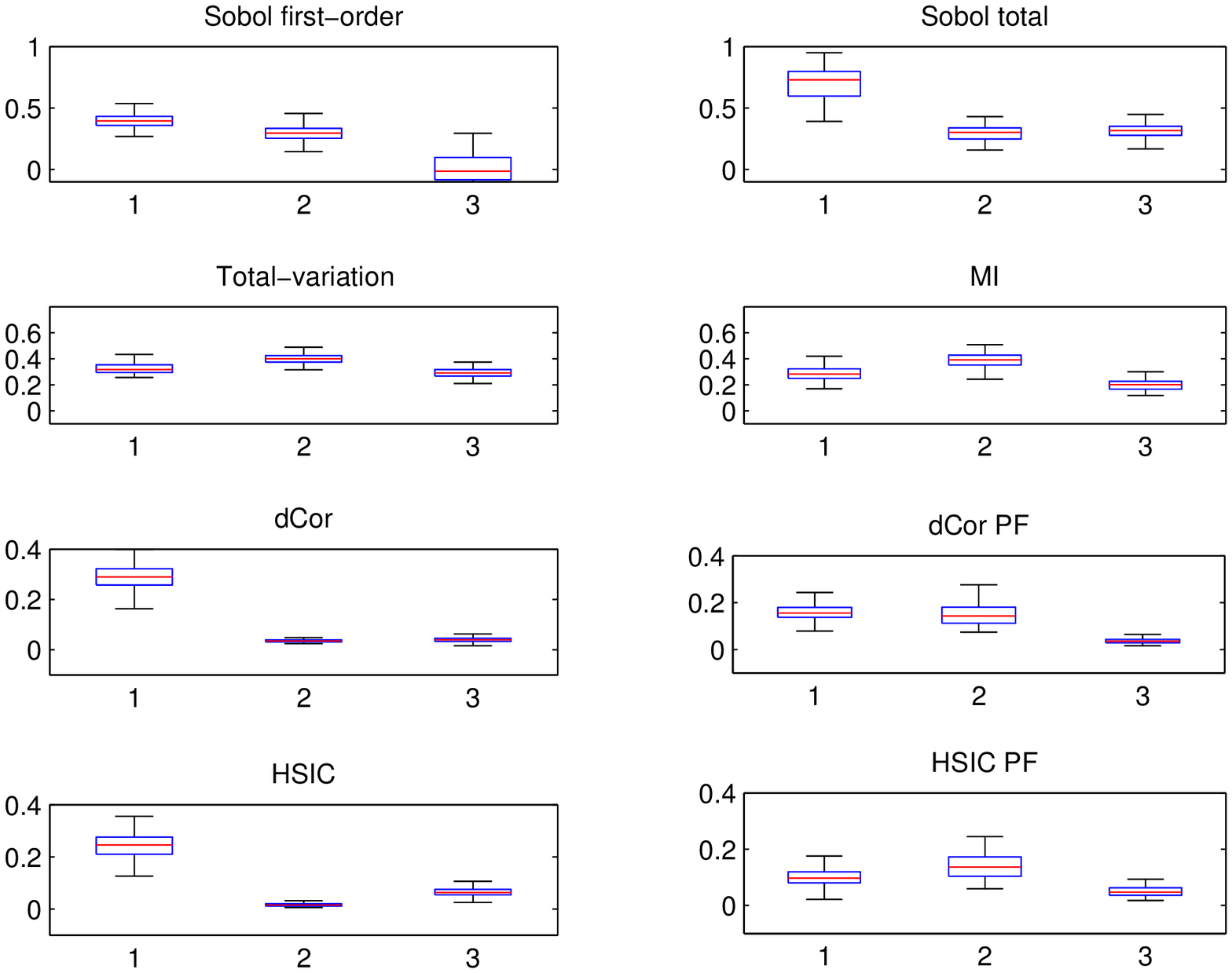}
\caption{First-order SI, total SI, $S_{X}^f$, $S_{X}^{dCor}$, $S_{X}^{dCorPF}$, $S_{X}^{HSIC_{\mathcal{F},\mathcal{G}}}$ and $S_{X}^{HSIC_{\mathcal{F},\mathcal{G}}PF}$ for function $\eta_3$, $n=200$, 100 replicates.} \label{figIshi1}
\end{figure}

This time we also compute $S_X^{TOT}$ since $\eta_3$ encompasses a strong iteration term. We use a sample of size $n=200$ for computing 
$S_X^1$, $S_{X}^f$, $S_{X}^{dCor}$ and $S_{X}^{dCorPF}$, but now we also need a sample of size $n\times p = 200\times 3=600$ for $S_X^{TOT}$ with RBD-FAST. Estimates obtained with 100 replications are reported in Figure \ref{figIshi1}. While first-order SIs indicate that $X_3$ has a negligible impact, it actually influences the output through an interaction term which is naturally accounted for by the total index. Is is interesting to note that all other indices detect the impact of $X_3$, as was pointed out by \citet{borgonovo07} for the total-variation index. However, one can observe the striking adequacy between $S_X^{TOT}$ and $S_{X}^{dCor}$ (unlike $S_X^f$). This  clearly shows that distance correlation has the potential to detect any interaction effect since it is specifically designed for nonlinear dependence. An additional appealing property is that its estimation does not depend on the number of inputs, unlike $S_X^{TOT}$. This is a major advantage for expensive computer codes. Finally, $S_{X}^{dCorPF}$ tends to bring the same information as $S_X^f$. But recall that it has the same limitation as $S_X^{TOT}$ concerning computational cost due to the pick-and-freeze technique. The same comments apply to $S_{X}^{HSIC_{\mathcal{F},\mathcal{G}}}$ and $S_{X}^{HSIC_{\mathcal{F},\mathcal{G}}PF}$.

\begin{figure}[h]
\centering
\includegraphics[width=0.9\textwidth]{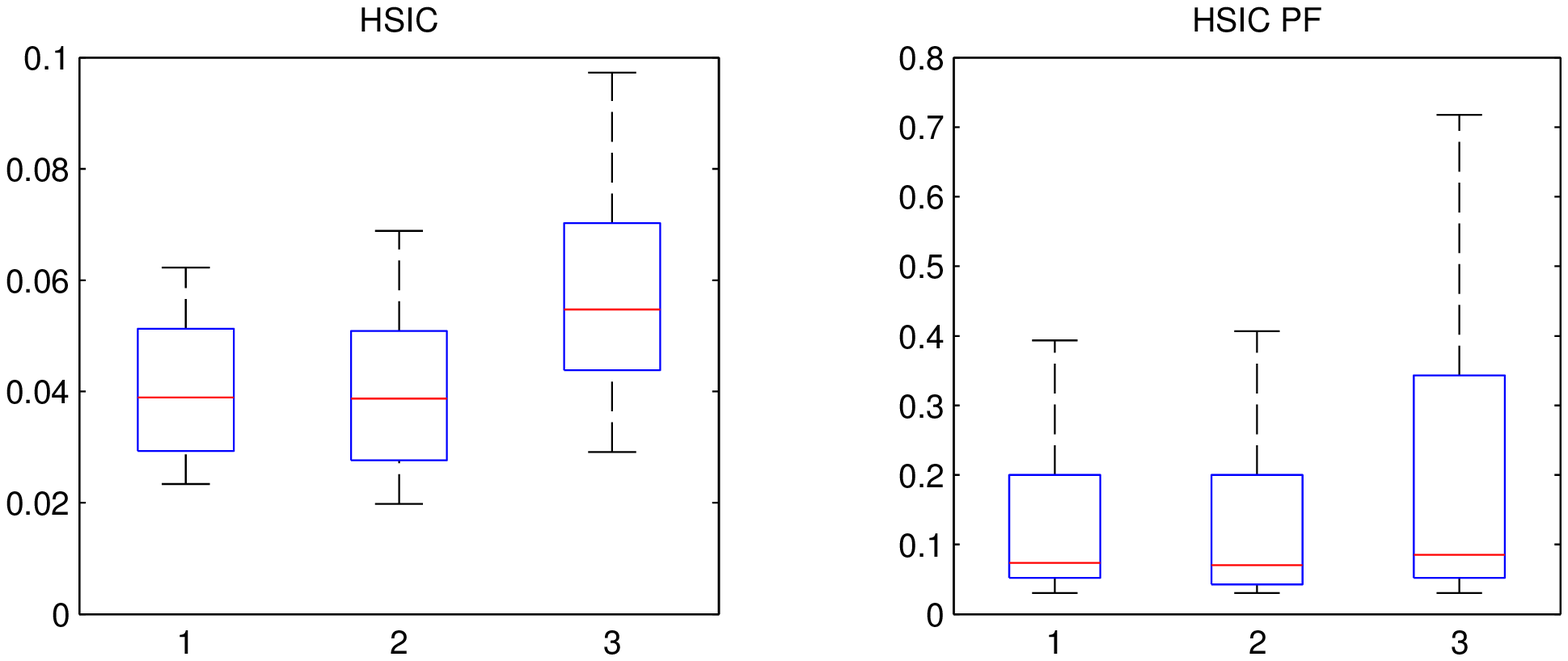}
\caption{$S_{X}^{HSIC_{\mathcal{F},\mathcal{G}}}$ and $S_{X}^{HSIC_{\mathcal{F},\mathcal{G}}PF}$ for function $\mathds{1}\{\eta_3>10\}$, $n=200$, 100 replicates.} \label{figIshi2}
\end{figure}

We also investigate HSIC on level sets of the Ishigami function to compare our results with \citet{fort13}. More precisely, we use a categorical kernel and use the change of variable $Z = \mathds{1}\{\eta_3(X)>10\}$. Figure \ref{figIshi2} shows that we can recover the fact that input factor $X_3$ is more important than $X_1$ and $X_2$ for this level set function, as was observed by \citet{fort13}.

\paragraph{Screening.} We now propose to study the performance of feature selection as an alternative to screening for problems where the number of input variables is large (more than $20$). We will deliberately limit the number of samples in order to be as close as possible to a real test case on an expensive code.
\begin{itemize}
\item \citet{morris06} function 
\[\eta_4(X) = \alpha \sum_{i=1}^k \left( X_i + \beta \prod_{i<j=2}^k X_iX_j\right)\] 
where $\alpha = \sqrt{12} -6\sqrt{0.1(k-1)}$, $\beta=\sqrt{12}\sqrt{0.1(k-1)}$, $X_i\sim \mathcal{U}(0,1)$, $i=1,\ldots,30$ and $1\leq k\leq 10$ is an integer controlling the number of influential inputs.
\end{itemize}

We select $n=50$, $k=5$ and compute $S_X^1$, $S_X^f$, $S_{X}^{dCor}$ and $S_{X}^{HSIC_{\mathcal{F},\mathcal{G}}}$ since all other indices are too expensive to compute in this setting (recall that $p=30$). First remark that first-order SIs identify the influential inputs in mean, but there are many replicates for which they are confounded with non-influential ones. On the contrary, $S_X^f$ completely fails ate detecting them: it will then be excluded from the other tests on screening. Notably, dCor and HSIC perfectly discriminate the first five factors and identify the remaining ones as independent from the output. We also use the HSIC Lasso (\ref{hsiclasso}), where for each replicate we use a bootstrap procedure to evaluate the probability of selection of each input factor. Here, HSIC Lasso performs very well since it selects the first 5 inputs factors almost every time.

\begin{figure}[h]
\centering
\includegraphics[width=0.9\textwidth]{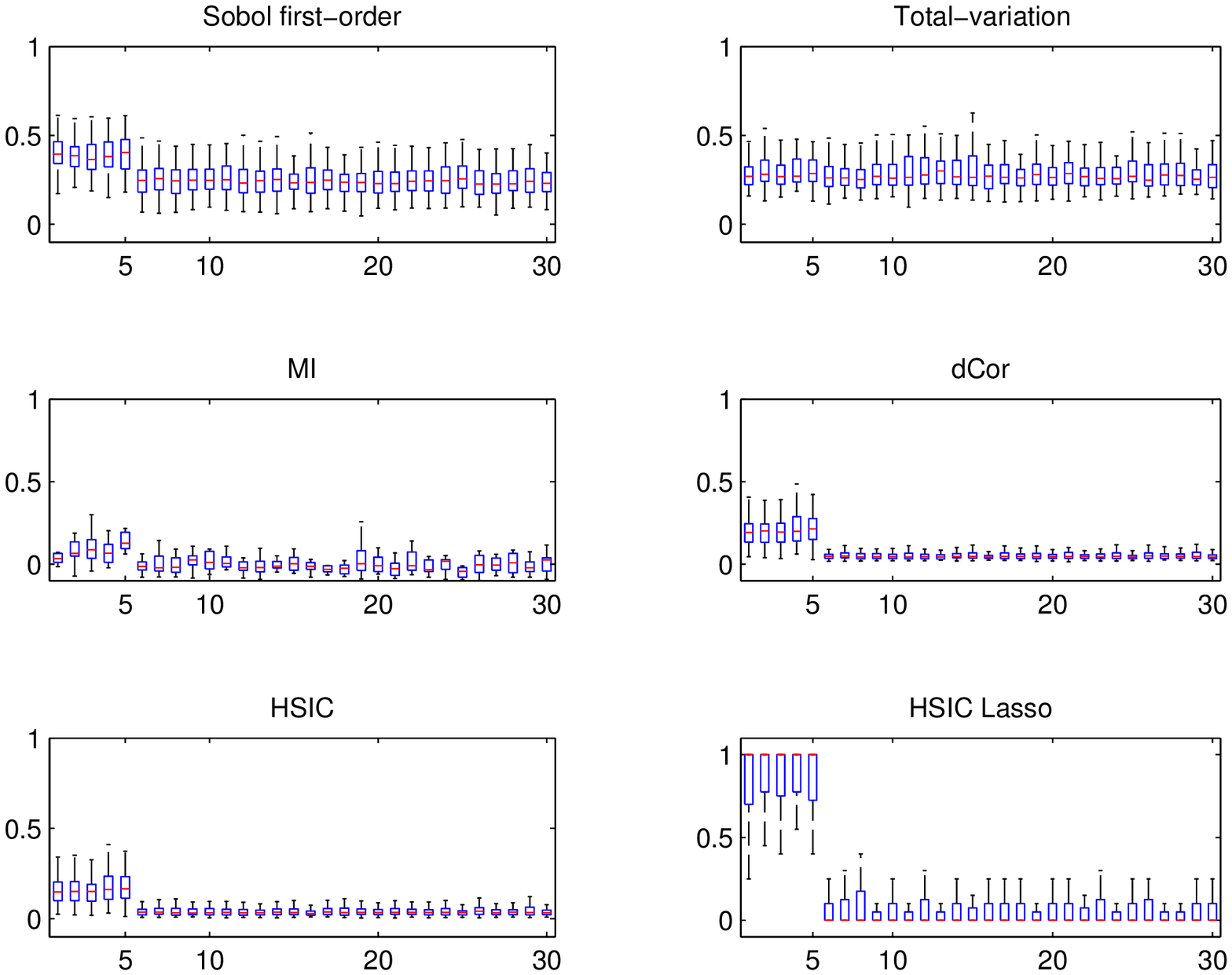}
\caption{First-order SI, $S_{X}^f$, $S_{X}^{dCor}$, $S_{X}^{HSIC_{\mathcal{F},\mathcal{G}}}$ and HSIC Lasso for function $\eta_4$, $n=50$, 100 replicates.} \label{figMorris1}
\end{figure}

\begin{itemize}
\item \citet{sobol99} function
\[\eta_5(X) =\exp\left(\sum_{i=1}^{20} b_iX_i\right) - \prod_{i=1}^p \frac{\exp(b_i)-1}{b_i}\]
with $X_i\sim \mathcal{U}(0,1)$, $i=1,\ldots,20$ and $b_i$ are taken from \citet{moon12}. Only the first eight factors are influential.
\end{itemize}

\begin{figure}[h]
\centering
\includegraphics[width=0.9\textwidth]{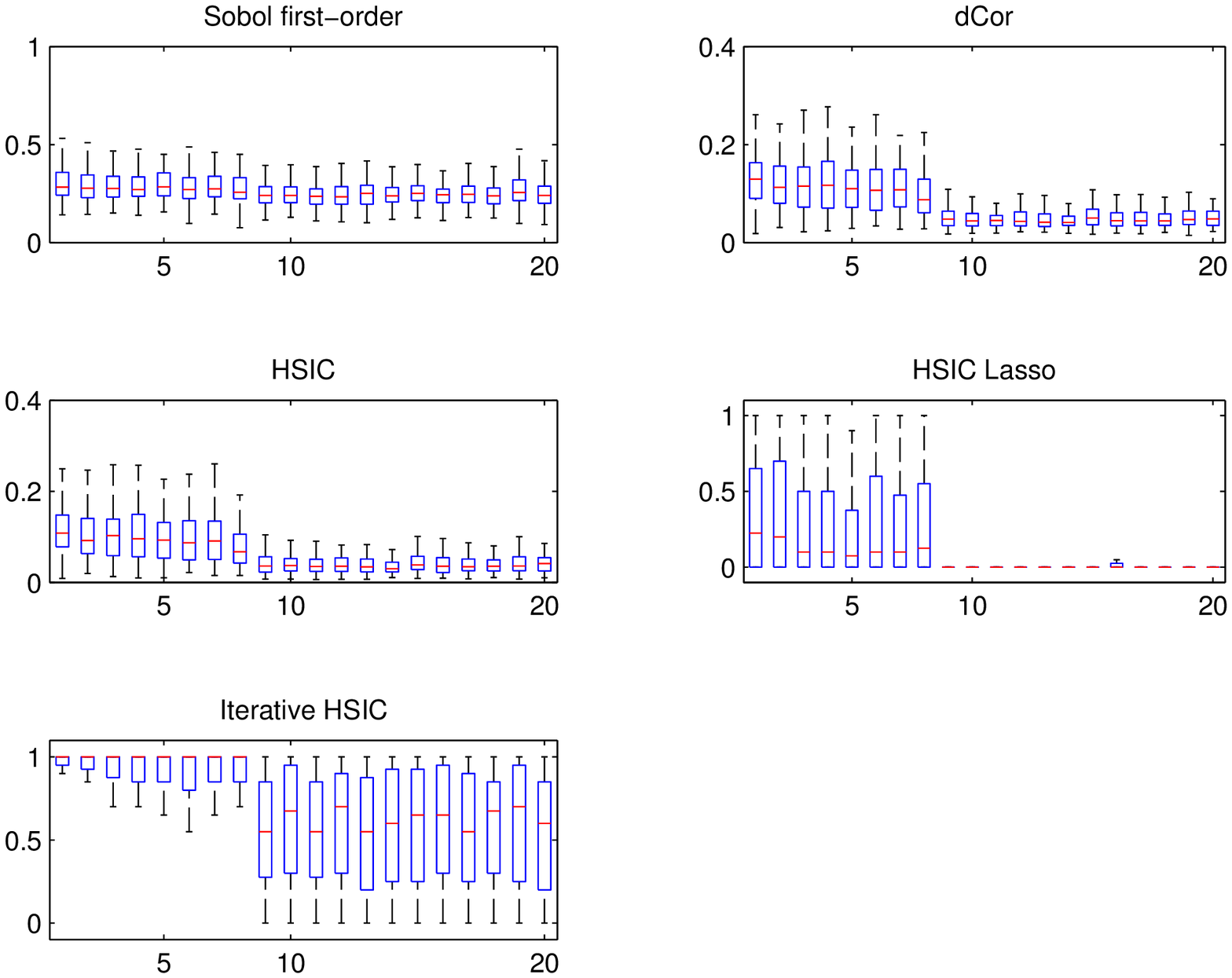}
\caption{First-order SI, $S_{X}^f$, $S_{X}^{dCor}$, $S_{X}^{HSIC_{\mathcal{F},\mathcal{G}}}$ and HSIC Lasso for function $\eta_5$, $n=50$, 100 replicates.} \label{figSoblev1}
\end{figure}

Here we also illustrate the feature selection method based on the iterative HSIC scheme detailed in Section \ref{secfs}.
Results for $n=50$ are given in Figure \ref{figSoblev1}. As expected, $S_X^1$ is unable to detect correctly the impact of the inputs. On the other hand, dCor and HSIC accurately estimate higher dependence for the first input factors than for the remaining ones. Similarly, HSIC Lasso never selects the last inputs as influential ones. The iterative feature selection based on HSIC performs well but tends to select more inputs than necessary. 

To go further and to compare our results with the ones obtained by \citet{moon12}, we repeat this experiment with $n=100$ (approximately the sample size used by \citet{moon12}). This time first-order SIs slightly detect the influential inputs, but the more interesting fact is that dCor, HSIC and HSIC Lasso give even better results and almost perfectly identifies them. Finally, the iterative HSIC scheme now almost always discards the non-influential inputs.

\begin{figure}[h]
\centering
\includegraphics[width=0.9\textwidth]{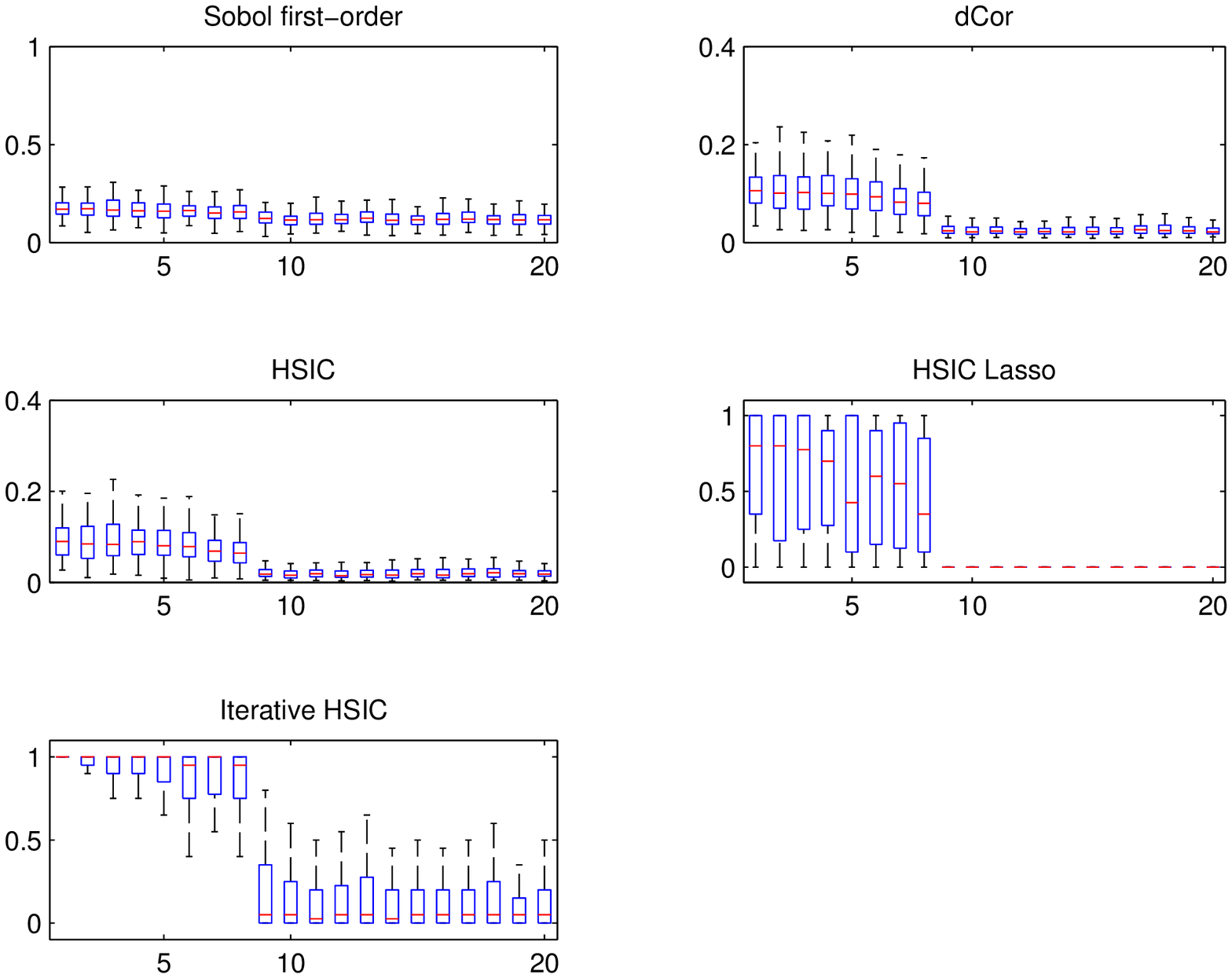}
\caption{First-order SI, $S_{X}^f$, $S_{X}^{dCor}$, $S_{X}^{HSIC_{\mathcal{F},\mathcal{G}}}$ and HSIC Lasso for function $\eta_5$, $n=100$, 100 replicates.} \label{figSoblev2}
\end{figure}

\subsection{Industrial applications}

\paragraph{Acquisition strategy for reservoir characterization.}
In the petroleum industry, reservoir characterization aims at reducing the uncertainty on some unknown physical parameters of an oil reservoir by using all the data collected on the field, e.g. well logs, seismic images or dynamic data at the wells (pressures, ...). Basically, engineers solve a Bayesian inverse problem where an initial prior distribution assumed on the parameters is updated by incorporating all field observations to produce a posterior distribution. In the end, this posterior distribution is used to predict the expected oil recovery of the reservoir in the future. 
However, it may be expensive to collect data and usually one wants to gather relevant observations only. This principle is at the core of so-called data acquisition strategies. For example, given the prior distribution, a natural idea is to get data which, when incorporated in the Bayesian procedure, will reduce the most the uncertainty of the obtained posterior distribution. It is easy to see that this idea actually corresponds to performing a sensitivity analysis of the parameters when data varies, where the difference between the prior and the posterior distribution is given by the measure of uncertainty reduction one chooses, i.e. the dissimilarity function  $d(\cdot,\cdot)$ in equation (\ref{Si}). Since the number of both uncertain parameters and observations can be large, we can greatly capitalize on the advantages of dCor and HSIC measures in arbitrary dimension to perform this task. 

\begin{figure}[h]
\centering
\includegraphics[width=0.9\textwidth]{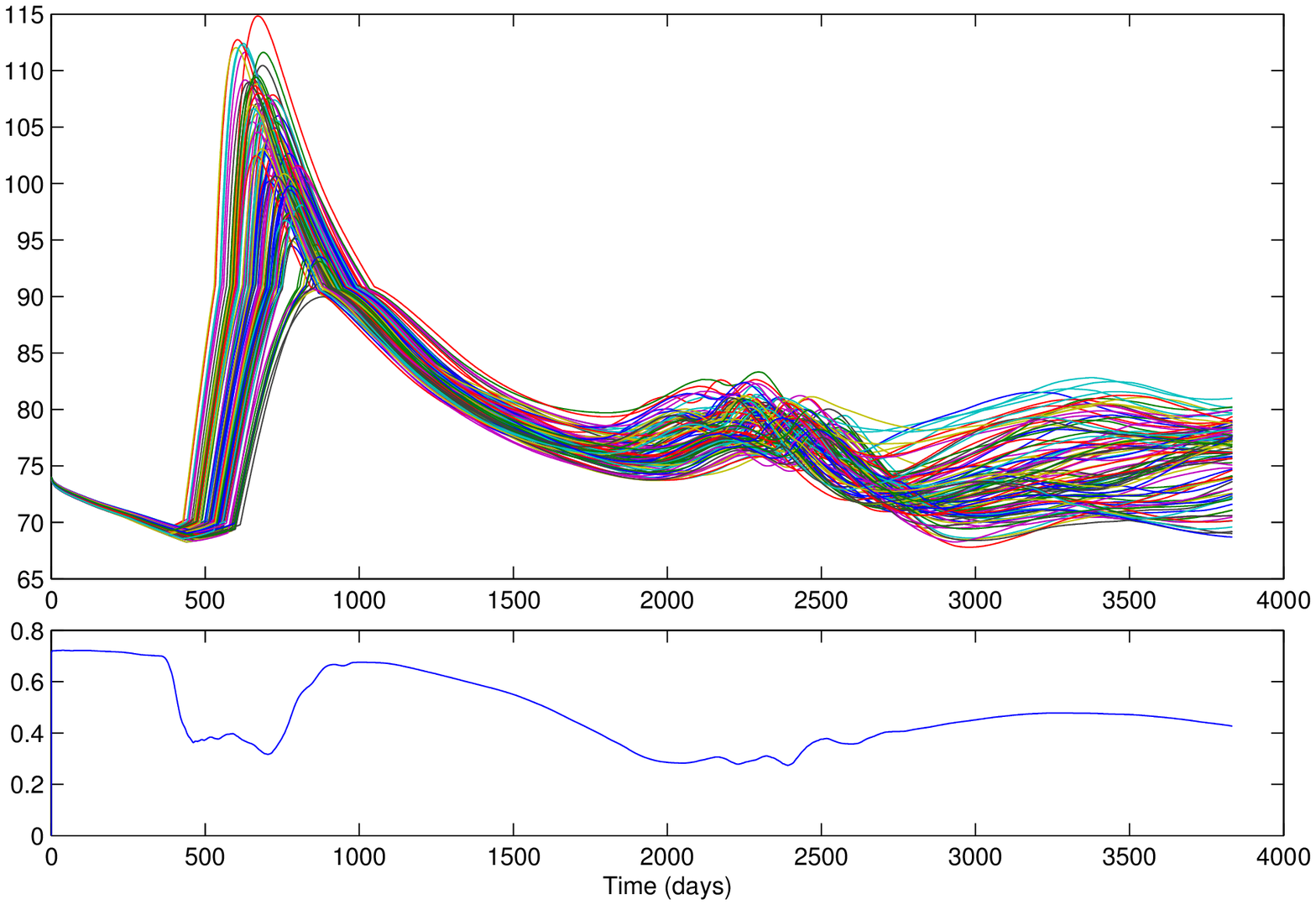}
\caption{Sample of size $n=100$ of the collectable data (top) and dCor measure between the parameters and the data at each time step (bottom) for the Punq test case.} \label{figpunq}
\end{figure}

Our example here makes use of the Punq reservoir test case, which is an oil reservoir model derived from real field data \citep{manceau01}. In this simplified model, seven variables which are characteristic of media, rocks, fluids or aquifer activity, are considered as uncertain (permeability multipliers, residual oil saturations, ...) and are assigned a uniform prior distribution. For illustration purposes, we assume that collectable data only consist of gas-oil ratios measurements at a given well. We generate a sample of size $n=100$ of the prior distribution, and propagate them through a fluid-flow simulator to get a sample of the simulated gas-oil ratios at the well over 10 years.  They are given in Figure \ref{figpunq}, top. For each day in these 10 years, we compute the dCor measure between the parameters and the simulated ratio, see Figure \ref{figpunq}, bottom. This information makes it possible to pick up the days where measurements should be collected in order to reduce as much as possible the uncertainty on the parameters: at the beginning of the reservoir production (before 500 days) or around 3 years after. Obviously, this procedure generalizes to any number of measurements and can be performed sequentially on many observations thanks to the properties of dCor.

\begin{figure}[h]
\centering
\includegraphics[width=0.9\textwidth]{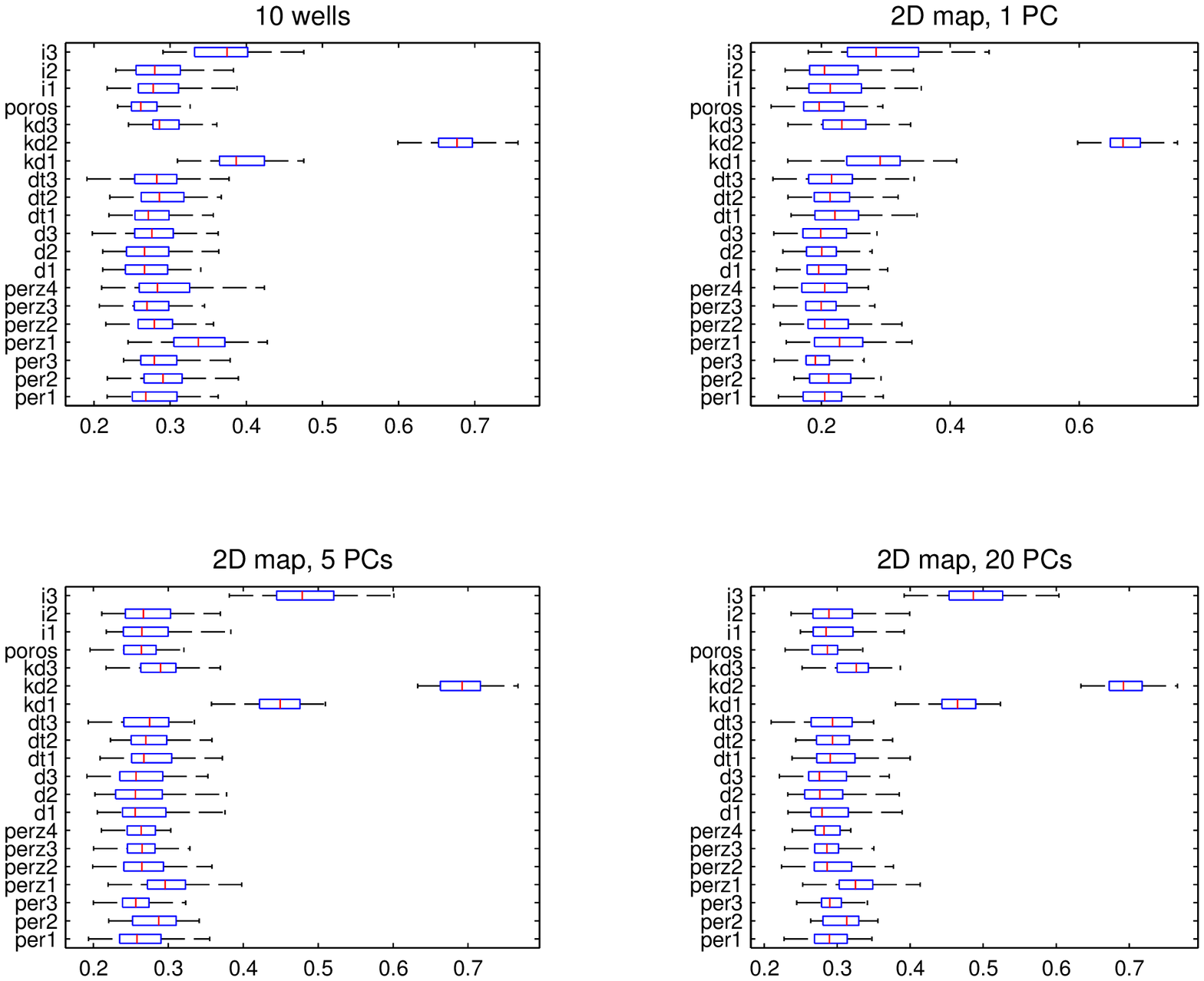}
\caption{HSIC measure between the parameters and the data for the Marthe test case (10 observation wells and 2D maps with a PCA kernel with 1, 5 and 20 principal components).} \label{figmarthe}
\end{figure}

\paragraph{Screening for contamination migration in waste storage site.}
The Marthe test case investigated here concerns prediction of the transport of strontium 90 in a porous water-saturated medium for evaluating the contamination of an aquifer in a temporary storage of radioactive waste \citep{volkova08}. Twenty input parameters mainly representative of the geological uncertainty are considered as random, and a set of 300 simulations is available at \url{http://www.gdr-mascotnum.fr/benchmarks.html}. Accessible outputs are strontium 90 concentrations simulated at ten different wells, as well as the concentrations on a complete 2D map of the area (discretized on $64\times 64 = 4096$ pixels). We place ourselves in a screening setting where we use only $n=50$ simulations to identify the influential inputs. To estimate the variability of our results, we pick at random these 50 samples among the 300 available and repeat the procedure 100 times. We use the HSIC dependence measure with a Gaussian kernel first in its standard form by considering the vector of concentrations at the 10 observation wells. But we also take advantage of a kernel designed for the 2D maps as was mentioned in Section \ref{secFunc}. Namely we use the PCA semi-metric \citep{ferraty06} and vary the number of principal components (1, 5 and 20 explaining $50\%$, $80\%$ and $95\%$ of the total variance, respectively). Results are given in Figure \ref{figmarthe}. First note that they are coherent with the ones obtained by \cite{volkova08} where the authors used the 300 simulations to build a surrogate model. Here, we then get the same detection of influential inputs but with only 50 simulations (parameters i3, kd1, kd2). In addition, the PCA kernel leads to more discriminating indices as soon as the explained variance is sufficient (5 PCs). This clearly illustrates the potential of HSIC for functional data.

\section{Conclusion}
In this paper, we introduced a new class of sensitivity indices based on dependence measures which overcomes the insufficiencies of variance-based methods in GSA. 
We demonstrated that when the output distribution is compared with its conditional counterpart through Csisz\'ar f-divergences, sensitivity indices arise as well-known dependence measures between random variables. We then extended these indices by using recent state-of-the-art dependence measures, such as distance correlation and the Hilbert-Schmidt independence criterion. We also emphasized the potential of feature selection techniques relying on such dependence measures as alternatives to screening in high dimension.

Interestingly, these new sensitivity indices are very robust to dimensionality, have low computational cost and can be elegantly extended to functional and categorical output or input variables. This opens the door to new and powerful tools for GSA and factors screening for high dimensional and expensive computer codes.

\bibliographystyle{agsm}

\end{document}